\documentclass[preprint,1p,11pt]{elsarticle}
\usepackage{mathrsfs}
\usepackage{epsfig}
\usepackage{amsmath}
\usepackage{amssymb}
\usepackage{slashed}

\newcommand{\bea}{\begin{eqnarray}}
\newcommand{\eea}{\end{eqnarray}}
\newcommand{\bnn}{\begin{eqnarray*}}
\newcommand{\enn}{\end{eqnarray*}}
\newcommand{\be}{\begin{equation}}
\newcommand{\ee}{\end{equation}}

\newcommand{\sn}{\,\mathrm{snc}\,}
\newcommand{\cn}{\,\mathrm{csc}\,}

\journal{Journal of Mathematical Analysis and Applications}
\def\MSC{\par\leavevmode\hbox {\it MSC:\ }}%

\begin{document}
\begin{frontmatter}

\title{Holographic coordinates}

%\author{S. Ulrych\\ Wehrenbachhalde 35, CH-8053 Z\"urich, Switzerland}
\author{S. Ulrych}
%\affiliation{Wehrenbachhalde 35, CH-8053 Z\"urich, Switzerland}
\address{Wehrenbachhalde 35, CH-8053 Z\"urich, Switzerland}
%\date{May 16, 2010}
%\maketitle

\begin{abstract}
The Laplace equation in the two-dimensional Euclidean plane is considered  in the context of the inverse stereographic projection.
The Lie algebra of the conformal group as the symmetry group of the Laplace equation can be represented solely in terms of the solutions and derivatives of the solutions of the Laplace equation. 
It is then possible to put contents from differential geometry and quantum systems, like the Hopf bundle, relativistic spin, bicomplex numbers, and the Fock space into a common context. 
The basis elements of the complex numbers, considered as a Clifford paravector algebra, are reinterpreted as differential tangent vectors referring to dilations and rotations.
In relation to this a homogeneous space is defined with the Lie algebra of the conformal group, where dilations and rotations are 
the coset representatives. Potential applications in physics are discussed.
\end{abstract}

%04.20.Cv Fundamental problems and general problems

\end{frontmatter}

\section{Introduction}
\label{intro}
The dimensional reduction of Minkowski space has been discussed more than two decades ago by  't Hooft \cite{Hoo93}.
These investigations lead to what is known today as the holographic principle \cite{Sus95}.
It has been suggested that the holographic principle is a foundation of a quantum gravity theory
in analogy to the equivalence principle as a foundation of general relativity \cite{Bou02}.
The holographic principle influenced string theory and lead to the
AdS/CFT correspondence \cite{Mal98,Gub98,Wit98}.
The holographic principle is applied in different physical areas like for example loop quantum gravity, cosmology, QCD, condensed matter physics, instantons, and
neutron stars \cite{Sak05,Sak05a,Bro06, Her06, Hat07, Her07, Har09, Sut10, Rho10, McG10, Ash11, San12, Kap15, Hoy16}. 

In accordance with the holographic principle 
one can see projective mappings as geometric operations, which transform
within representation spaces of different dimension.
Here M\"obius geometries and projective Lorentz spaces play a central role. For the mathematical 
background with respect to M\"obius geometries it is referred to Hertrich-Jeromin \cite{Her03}, Sharpe \cite{Sha97}, Kisil \cite{Kis12,Kis16}, or Jensen et al. \cite{Jen16}.
Furthermore, a hierarchy of projective geometric spaces and M\"obius geometries
can be introduced based on Vahlen matrices \cite{Vah02,Ahl85,Ahl86, McI16}.
Maks \cite{Mak89} investigated explicitly the sequence
of M\"obius geometries represented in terms of Clifford algebras, which is considered also in \cite{Ulr17} in a bicomplex matrix representation.
In relation to particle physics one may assume that projections of a more general geometry can arise as part of the measurement process via electromagnetic forces.
In accordance with the holographic principle M\"obius geometries and M\"obius transformations can provide here a new perspective on the experimental data.

Atiyah, Manton, and Schroers \cite{Ati12, Ati16, Ati17, Ati18} describe electrons, protons, neutrinos, and neutrons within a model that has been inspired by Skyrme's baryon theory \cite{Sky61}.
The Skyrme model became popular when Witten came up with the idea that baryons arise as solitons of the classical meson fields \cite{Wit79}.
The Skyrme theory turned out to be remarkably successful in describing baryons in the previous decades \cite{Rho16a}.
Thus matter is supposed to arise from a pure geometric foundation \cite{Man04, Dun10}. 
Solitons are usually considered in the context of gauge transformations, which relate in a mathematical context to fiber bundles and Cartan geometries \cite{Sha97}.

With regards to the topics discussed above research in two-dimensional conformal field theories \cite{Bel84,Gin88,Blu09, Kis07,Mur17} has turned out to be important,
also because conceptual insights can be obtained more easily in a lower dimensional geometry.
Motivation to proceed in this direction is also provided by \cite{Ulr17}. 
In this article the two-dimensional Euclidean plane and the conformal symmetries, which refer to the Laplace equation in $R^2$, are fundamental in a
sequence of higher dimensional geometries and Clifford algebras.
Therefore, one can ask the question how the Laplace equation can be embedded into higher dimensional geometries.
Here concepts like the stereographic projection should play an important role.
In this sense one can consider the Laplace equation in the base space $S^2$ instead of $R^2$ by means of what is denoted here as holographic coordinates. 
This conceptual change is adopted from Lie sphere geometry \cite{Jen16,Cec92}, where the spherical geometry is leading to considerable simplifications compared to the Euclidean geometry. 

The main result of this article is that with the help of the holographic coordinates
the Lie algebra of the conformal group in two dimensions can be represented solely in terms of the solutions and derivatives of the solutions of the corresponding Laplace equation. 
It turns out that different concepts in mathematics and physics, like projective spaces, spin, the Hopf bundle, Minkowski space, or bicomplex numbers can be seen in relation to each other.
The connections are provided by the solutions and symmetries of the Laplace equation.

\section{Conformal transformations}
\label{Conftrafo}
Conformal transformations have been investigated in \cite{Ulr17} based on spin representations over the ring of bicomplex numbers.
The non-zero commutation relations for these spin matrices are
\bea
\label{confcomm}
\left[s_{\mu\nu},p_\sigma\right]&=&g_{\nu \sigma}p_{\mu}-g_{\mu \sigma}p_{\nu}\,,\nonumber\\
\left[s_{\mu\nu},q_\sigma\right]&=&g_{\nu \sigma}q_{\mu}-g_{\mu \sigma}q_{\nu}\,,\nonumber\\
\left[b,p_\mu\right]&=&-p_\mu\,,\nonumber\\
\left[b,q_\mu\right]&=&q_\mu\,,\nonumber\\
\left[q_\mu,p_\nu\right]&=&2(g_{\mu\nu}b+s_{\mu\nu})\,.
\eea
Here $s_{\mu\nu}$ denotes the angular momentum operator, $p_\mu$ the momentum,
and $q_\mu$ the special conformal transformations.
Compared to \cite{Ulr17} the symbol $b$ instead of $d$ is used to label the scale transformation in order to avoid future confusion with the exterior derivative.
The commutation relation between two angular momentum operators is not listed here. It will be discussed separately later in the text.

The difference of the matrix commutation relations compared to the relations used by Kastrup \cite{Kas62} is that no
complex units appear on the right-hand side of the equations. If one tries to represent the spin matrices
in terms of differential operators, one cannot apply the conventions of quantum physics and add complex units
to the differential operators. However, one can start in the context of classical field theory and use
the conventions of Freedman and Van Proeyen \cite{Fre12}, which result in the following representation
\bea
\label{diffrep}
s_{\mu\nu}&=&x_\mu\partial_\nu-x_\nu\partial_\mu\,,\nonumber\\
p_\mu&=&\partial_\mu\,,\nonumber\\
q_\mu&=&2x_\mu x_\nu \partial^\nu-x^2\partial_\mu\,,\nonumber\\
b&=&x_\mu\partial^\mu\,.
\eea
The differential operators refer to the symmetries of the Laplace equation and they thus all commute with
the Laplace operator.

\section{Laplace operator}
\label{Laplace}
One can reconsider now the method discussed in \cite{Ulr14} in the simplified context of the two-dimensional plane.
Nevertheless, a notation is used which applies formally in arbitrary metric spaces.
The Laplace operator is written as
\be
\label{simplaplace}
\triangle=p\cdot p\,.
\ee
The momentum operator is expanded in terms of the Clifford  basis vectors $p=p^\mu e_\mu$.
The dot product is evaluated with the help of the metric tensor, which is defined by
\be
\label{metric}
g_{\mu\nu}=\frac{e_\mu\bar{e}_\nu+e_\nu\bar{e}_\mu}{2}\,.
\ee
As discussed in the previous section the imaginary units will be removed from the momentum operator in contrast to the representation in \cite{Ulr14}
\be
p^\mu=\partial^\mu=\partial_{x_\mu}\,.
\ee
In the context of the Laplace equation one may simplify the Laplace operator to $p\bar{p}$.
It is supposed that the momentum operator is invariant with respect to the Clifford conjugation $\bar{p}^\mu=p^\mu$.
Finally, in the simple two dimensional case the paravector model $e_\mu=(1,e_1)$ is used
with the single non-trivial basis element $e_1=i$.
One arrives at
\be
\triangle=\left(\partial_{x_0}-i \partial_{x_1}\right)
\left(\partial_{x_0}+i \partial_{x_1}\right).
\ee
This can be further evaluated to obtain the standard representation of the Laplace operator in Cartesian coordinates.
The sign conventions for $p$ and $\bar{p}$ are aligned with \cite{Gue08}.

\section{Polar coordinates}
\label{polarrep}
The polar coordinates stand in the center of many applications within physics.
This section repeats well known properties of the polar coordinates.
The intention is to let the reader become familiar with the notation used in this article.
This allows also to compare with
the changes that arise when the holographic coordinates are introduced in the following sections.

The coordinate functions for the polar coordinates are introduced by the following definition
\be
\label{para}
x_{\mu}=\left(\begin{array}{c}
x_{0}\\
x_{1}\\
\end{array}\right)=
\left(\begin{array}{c}
r\cos{\varphi}\\
r\sin{\varphi}\\
\end{array}\right).
\ee
Here $r$ corresponds to the radius and $\varphi$ to the usual angle parameter in the two-dimensional plane \cite{Arn89}.
The new coordinates can be placed in the vector
\be
\label{sphcoord}
y_{\alpha}=\left(\begin{array}{c}
y_{0}\\
y_{1}\\
\end{array}\right)=
\left(\begin{array}{c}
r\\
\varphi\\
\end{array}\right).
\ee
In the following the indices $\alpha$ and $\beta$ are used in curved space.
The indices $\mu$ and $\nu$ are applied in flat space to indicate that the coordinates refer to another metric.

Basis vectors in the polar coordinate system can be introduced with the following definition
\be
\label{basedef}
\varsigma^\alpha=\frac{\partial x}{\partial y_\alpha}\;.
\ee
With the coordinate functions introduced above 
the new basis vectors can be calculated explicitly
\be
\label{sphbasis}
\varsigma^{0}=
\left(\begin{array}{c}
\cos{\varphi}\\
\sin{\varphi}\\
\end{array}\right),\quad
\varsigma^{1}=
\left(\begin{array}{c}
-r\sin{\varphi}\\
r\cos{\varphi}\\
\end{array}\right).
\ee
The result is used to calculate the metric tensor which is given by the scalar product of
the basis vectors. One finds
\be
g^{\alpha\beta}= \varsigma^\alpha \cdot \varsigma^\beta=\left(\begin{array}{cc}
1&0\\
0&r^2\\
\end{array}\right).
\ee
The metric tensor can be applied to raise the index of the above coordinate vector. 

The transformation matrix referring to a coordinate change can be calculated as
\be
\label{trafomat}
A_\mu^{\phantom{\mu}\alpha}=\frac{\partial x_\mu }{\partial y_\alpha}=\left(\begin{array}{cc}
\cos{\varphi}&-r\sin{\varphi}\\
\sin{\varphi}&r\cos{\varphi}\\
\end{array}\right).
\ee
One can raise and lower indices with the corresponding metric tensors.
The first metric tensor is referring to flat space, the second metric tensor to curved space
\be
\label{cindextrafo}
A^\mu_{\phantom{\mu}\alpha}=g^{\mu\nu}g_{\alpha\beta}A_\nu^{\phantom{\nu}\beta}\,.
\ee
The metric tensor $g^{\mu\nu}$ is simply the identity matrix as it refers to the Euclidean plane.
The second metric tensor corresponds to the inverse of $g^{\alpha\beta}$.
The pieces can be put together and one obtains the transformation matrix
\be
\label{lamtrafo}
A^\mu_{\phantom{\mu}\alpha}=\left(\begin{array}{cc}
\cos{\varphi}&-r^{-1}\sin{\varphi}\\
\sin{\varphi}&r^{-1}\cos{\varphi}\\
\end{array}\right).
\ee
This completes the brief overview of the polar coordinates and the notation used in this article. 

\section{Laplace equation}
\label{Laplaceeq}
The polar coordinates introduced in the previous section can now be applied to transform the Laplace operator.
The Laplacian in polar coordinates can be calculated as usual, nevertheless a redefined form is useful in the following
\be
\label{pollap}
\triangle=r\partial_rr\partial_r
+\partial_{\varphi\varphi}\,.
\ee
To recover the standard representation of the Laplace operator one has to scale from the left side by $r^{-2}$.
The rescaling can be applied in the context of the scale invariant Laplace equation.

From the physical point of view using the Laplace equation as a starting point is understood in the sense of a zero-mass pre-gauge theory in an Euclidean worldsheet
\be
\triangle\upsilon=0\,.
\ee
Masses are supposed to arise at a later stage through symmetry breaking of the chiral symmetry for the strong force and through the Higgs mechanism 
for the weak interaction \cite{Gui91,Wil12, Rho16}.
The solution of the Laplace equation with separation constant $\alpha$
is given by
\be
\upsilon^\alpha=e^{i\alpha\varphi}r^{\alpha}\,.
\ee
It should be clear from the context when $\alpha$ denotes a vector index or the separation constant of the Laplace equation.
The separation constant can also be complex valued. The solutions for $\alpha=0$ can be found in the literature.

\section{Conformal compactification}
%\label{compact}
In order to motivate the holographic coordinates, it is worth to look at the conformal compactification of the two-dimensional Euclidean space.
In addition to the coordinates $x_0$ and $x_1$, which are copied to the new conformal vector $u_\mu$, two further
coordinates are introduced within the two-dimensional geometry of the base space.
The new coordinates correspond to a rescaling of the inverse stereographic projection from the South pole represented within the projective space $\mathbb{RP}^3$, see for example Hertrich-Jeromin \cite{Her03} or Cecil~\cite{Cec92}
\be
u_\mu=\left(\begin{array}{c}
x_0\\
x_1\\
(1-x^2)/2\\
(1+x^2)/2
\end{array}\right)=\left(\begin{array}{c}
2x_0/(1+x^2)\\
2x_1/(1+x^2)\\
(1-x^2)/(1+x^2)\\
1
\end{array}\right).
\ee
With the help of this representation one can motivate the introduction of the holographic coordinates,
which will be discussed in the following sections. The vector on the right-hand side can now be reparametrized as
\be
u_\mu=
\left(\begin{array}{c}
\sin{\theta}\cos{\varphi}\\
\sin{\theta}\sin{\varphi}\\
\cos{\theta}\\
1
\end{array}\right).
\ee
The vectors are situated in Minkowski space $\mathbb{R}^{3,1}$ and have all zero length.
One thus finds a correspondence between the original coordinates in the plane and coordinates on the sphere $S^2$.

\section{Holographic coordinates}
\label{holcord}
The two-dimensional plane can again be projected out of the four-dimensional vectors discussed in the previous section.
The vector in the plane is thus represented in terms of the following coordinates
\be
\label{parahol}
x_{\mu}=\left(\begin{array}{c}
x_{0}\\
x_{1}\\
\end{array}\right)=
\left(\begin{array}{c}
\sin{\theta}\cos{\varphi}\\
\sin{\theta}\sin{\varphi}\\
\end{array}\right).
\ee
The parameter $\varphi$ remains here the same as in the case of polar coordinates.
The radius $r$ will be replaced by an angle parameter $\theta$.
The angle parameter $\theta$
refers to the intersection point on the two-dimensional sphere $S^2$, which is introduced as part of the inverse stereographic projection,
with the projection line to the coordinate vector in the two-dimensional plane $\mathbb{R}^2$.
Consider here also Figure 7.2 in Freedman and Van Proeyen \cite{Fre12}.

The new coordinate vector is given by
\be
%\label{sphcoord}
y_{\alpha}=\left(\begin{array}{c}
y_{0}\\
y_{1}\\
\end{array}\right)=
\left(\begin{array}{c}
\theta\\
\varphi\\
\end{array}\right).
\ee
One may apply now the formulas of Section \ref{polarrep} in order to calculate the metric tensor and the coordinate transformation matrices.
The metric tensor is calculated as
\be
g^{\alpha\beta}=\left(\begin{array}{cc}
\cos^2{\theta}&0\\
0&\sin^2{\theta}\\
\end{array}\right).
\ee
The transformation matrix of Eq.~(\ref{trafomat}) is recalculated in the new coordinate system again by
direct derivation of the coordinate functions with respect to the new coordinates
\be
A_\mu^{\phantom{\mu}\alpha}=\left(\begin{array}{cc}
\cos{\theta}\cos{\varphi}&-\sin{\theta}\sin{\varphi}\\
\cos{\theta}\sin{\varphi}&\sin{\theta}\cos{\varphi}\\
\end{array}\right).
\ee
The transformation matrix with raised and lowered indices
is calculated as discussed in Section \ref{polarrep} with the help of the inverse metric tensor $g_{\alpha\beta}$.
Again one can apply Eq.~(\ref{cindextrafo}).
This results in the following transformation matrix
\be
\label{holtrafo}
A^\mu_{\phantom{\mu}\alpha}=\left(\begin{array}{cc}
\cos^{-1}{\theta}\cos{\varphi}&-\sin^{-1}{\theta}\sin{\varphi}\\
\cos^{-1}{\theta}\sin{\varphi}&\sin^{-1}{\theta}\cos{\varphi}\\
\end{array}\right)\,.
\ee
The Laplace operator in two-dimensional Euclidean space then takes a simple form,
which includes a rescaling as discussed in Section~\ref{Laplaceeq}
\be
\label{hollap}
\triangle=\tan{\theta}\partial_\theta\tan{\theta}\partial_\theta
+\partial_{\varphi\varphi}\,.
\ee
The rescaling can be done in the context of the Laplace equation, because the Laplace equation is
invariant with respect to conformal transformations.

\section{Holographic solutions of the Laplace equation}
The solutions of the Laplace equation in holographic coordinates can be determined easily.
One solution is given for example by
\be
\label{SolofLap}
\upsilon^\alpha=e^{i\alpha\varphi}\sin^\alpha{\theta}\,.
\ee
The solution can be transformed with the conformal transformations into further solutions of the Laplace equation.

Due to the choice of the coordinate system it is not surprising that the solutions for non-negative integers $\alpha=l$ correspond to a subset of the spherical harmonics $Y^m_{l}(\theta,\varphi)$.
One finds that the solutions of the Laplace equation in holographic coordinates are proportional to
\be
\label{spharm}
\upsilon^l\sim Y^l_l(\theta,\varphi)\,.
\ee
Thus the angular and magnetic quantum numbers of the spherical harmonics are identical. Furthermore, the 
spherical harmonics with $m=-l$ are solutions of the Laplace equation in holographic coordinates.

\section{Conformal group}
\label{geointeract}
One can investigate how the Lie algebra of the conformal group acts
on the solutions of the Laplace equation. The elements of the Lie algebra as defined in Eq.~(\ref{diffrep}) can be represented in terms of
the holographic coordinates
\bea
\label{allop}
b&=&\tan{\theta}\,\partial_\theta\,,\nonumber\\
s_{01}&=&\partial_{\varphi}\,,\nonumber\\
p_0&=&\cos{\varphi}\cos^{-1}{\theta}\partial_\theta -\sin{\varphi}\sin^{-1}{\theta}  \partial_\varphi\,,\nonumber\\
p_1&=&\sin{\varphi}\cos^{-1}{\theta}\partial_\theta+\cos{\varphi}\sin^{-1}{\theta} \partial_\varphi\,,\nonumber\\
q_0&=&\cos{\varphi}\sin{\theta}\tan{\theta}\,\partial_\theta +\sin{\varphi}\sin{\theta}\, \partial_\varphi\,,\nonumber\\
q_1&=&\sin{\varphi}\sin{\theta}\tan{\theta}\,\partial_\theta-\cos{\varphi}\sin{\theta}\, \partial_\varphi\,.
\eea
Here the transformation matrix calculated in Eq.~(\ref{holtrafo}) has been used, for example $b=x_\mu A^\mu_{\phantom{\mu}\alpha}\partial^\alpha$ with $x_\mu$ given by Eq.~(\ref{parahol}).
The Lie algebra of the conformal group can act on the solutions of the Laplace equation as given by Eq.~(\ref{SolofLap}), which results in
\bea
\label{repchange}
b\upsilon^\alpha&=&\alpha\upsilon^\alpha\,,\nonumber\\
s_{01}\upsilon^\alpha&=&i\alpha\upsilon^\alpha\,,\nonumber\\
p_0\upsilon^\alpha&=&\alpha\upsilon^{\alpha-1}\,,\nonumber\\
p_1\upsilon^\alpha&=&i\alpha\upsilon^{\alpha-1}\,,\nonumber\\
q_0\upsilon^\alpha&=&\alpha\upsilon^{\alpha+1}\,,\nonumber\\
q_1\upsilon^\alpha&=&-i\alpha\upsilon^{\alpha+1}\,.
\eea 
The solutions of the Laplace equation in holographic coordinates are eigenfunctions of the generators of rotation and scale transformation.
Note that the transformations in Euclidean geometry, given by a rotation and a translation, can be represented alternatively by a rotation and a scale transformation.

In conformal field theory $\alpha$ is denoted as the scale dimension.
One finds here that the momentum increases the scale dimension, whereas special conformal transformations reduce the scale dimension by one unit,
see for example Sleight \cite{Sle16}. Here it is the other way around.

It is possible to combine momentum and special conformal transformations into Clifford paravectors $p$ and $\bar{q}$ and write the above
equations in a more compact form. Because the $\upsilon^\alpha$ are holomorphic functions $\bar{p}$ and $q$ annihilate the solutions of the Laplace equation.

\section{Special conformal transformations}
The special conformal transformations appear  in Eq.~(\ref{repchange}) as the counterpart of the translations.
It is therefore worth to note that the special conformal transformations can be interpreted as translations with regards to infinity as the origin.
This interpretation is possible because the special conformal transformations can be decomposed into inversion, translation, followed by reinversion \cite{Olv93}.
The inversion maps infinity to zero and vice versa.

\section{Minkowski space}
The conformal transformations acting on the solutions of the Laplace equation
can be related to rotation operators in a higher dimensional space with the following
relations \cite{Kas62}
\bea
\label{coprot}
p_\mu&=&-s_{\mu 2}-s_{\mu 3}\,,\nonumber\\
q_\mu&=&s_{\mu 2}-s_{\mu 3}\,,\nonumber\\
b&=&s_{23}\,.
\eea
The commutation relations given in Eq.~(\ref{confcomm}) can now be summarized in one single equation
\be
\label{comm}
[s_{\mu\nu},s_{\rho\sigma}]=g_{\mu\sigma}s_{\nu\rho}-g_{\mu\rho}s_{\nu\sigma}-g_{\nu\sigma}s_{\mu\rho}
+g_{\nu\rho}s_{\mu\sigma}\,.
\ee
The metric tensor in the resulting four-dimensional space can then be derived from operators,
which have been originally defined in the two-dimensional Euclidean space $\mathbb{R}^{2}$.
The metric induced in this way corresponds to the metric of Minkowski space~$\mathbb{R}^{3,1}$
\be
\label{relmetric}
g_{\mu\nu}=\left(\begin{array}{cccc}
1&0&0&0\\
0&1&0&0\\
0&0&1&0\\
0&0&0&-1\\
\end{array}\right).
\ee
What happens behind the scenes is that a model for the M\"obius plane in two-dimensional space has been set up with the above definitions.
The Lorentz group $SO(3,1,\mathbb{R})$, which is represented by the operators $s_{\mu\nu}$, is used to perform transformations within the M\"obius plane \cite{Sha97}.

\section{Generalized coordinate space created by solutions of a classical field equation}
Prolongation is a standard concept in the theory of differential equations \cite{Olv93,Ste90}.
Here the solution of the differential equation and its derivatives appear in the definition of the symmetry group of the differential equation.
In fact the Lie algebra of the conformal group as the symmetry group of the Laplace equation can be fully represented in terms of its solution and its derivatives alone.
From Eq.~(\ref{repchange}) one can derive the following relations
\bea
\label{difform}
b&=&\upsilon\partial_\upsilon\,,\nonumber\\
s_{01}&=&i\upsilon\partial_\upsilon\,,\nonumber\\
p_0&=&\partial_\upsilon\,,\nonumber\\
p_1&=&i\partial_\upsilon\,,\nonumber\\
q_0&=&\upsilon^2\partial_\upsilon\,,\nonumber\\
q_1&=&-i\upsilon^2\partial_\upsilon\,.
\eea 
The operators $b$, $p_0$, and $q_0$ can be reinterpreted as acting on the real projective line $\mathbb{RP}^1$ made up of the coordinate $\upsilon$.
According to Olver \cite{Olv93} the symmetry generators can be identified with the matrices
\be
\label{realproj}
b=\frac{1}{2}\left(\begin{array}{cc}
1&0\\
0&-1
\end{array}\right),\quad
p_0=\left(\begin{array}{cc}
0&1\\
0&0
\end{array}\right),\quad
q_0=\left(\begin{array}{cc}
0&0\\
-1&0
\end{array}\right).
\ee
The matrices correspond to the infinitesimal generators of the special linear group $SL(2,\mathbb{R})$.
The matrix representation of $SL(2,\mathbb{R})$ acts on the projective line represented in terms of the vector
\be
\upsilon\equiv
\left(\begin{array}{c}
\upsilon\\
1\\
\end{array}\right).
\ee
The correspondence of the matrix representations to linear fractional transformations is then implied as usual by the properties of the projective space
\be
\upsilon\mapsto \frac{\alpha\upsilon+\beta}{\gamma\upsilon+\delta},\quad 
\left(\begin{array}{cc}
\alpha&\beta\\
\gamma&\delta
\end{array}\right)
\in SL(2,\mathbb{R})\,.
\ee
For example one can calculate how the symmetry transformation generated by the matrix $q_0$ acts infinitesimally on the solution of the Laplace equation
\be
\label{qpole}
\upsilon+\varepsilon q_0\upsilon=\frac{\upsilon}{1-\varepsilon\upsilon}\approx \upsilon +\varepsilon \upsilon^2\,.
\ee
As expected one obtains the same result as if one operates with the vector field $q_0=\upsilon^2\partial_\upsilon$. 
The important point here is that the coordinate is given by the solution of the Laplace equation.
Thus one finds that the solutions of a classical field equation suddenly appear as coordinate functions $\upsilon$ in a generalized coordinate space.

\section{Spin in Minkowski space}
The operators $s_{01}$, $p_1$, and $q_1$ complete the picture.
They complexify the real projective representation, which then leads to the complex projective line $\mathbb{CP}^1$
and the complex group of linear fractional transformations $SL(2,\mathbb{C})$. The group $SL(2,\mathbb{C})$
corresponds to the spin group of the Lorentz transformations in Minkowski space. The explicit spin matrices can be constructed by means of
Eqs.~(\ref{coprot}) and (\ref{realproj}) together with the following relations, which are deduced from Eq.~(\ref{difform}) 
\be
\label{complexline}
s_{01}=ib\,,\quad p_1=ip_0\,,\quad q_1=-iq_0\,.
\ee
The special linear group can be decomposed as usual into two unitary groups $SL(2,\mathbb{C})=SU(2,\mathbb{C})\otimes SU(2,\mathbb{C})$.
This shows how a property of quantum systems like spin emerges from a classical field equation. In addition, the complex projective line $\mathbb{CP}^1$
is referring to the invariance with regards to phase transformations, a further feature of quantum systems.

\section{Higher spin}
Higher spin representations have been considered in recent years in the context of string theory, investigations of the matter spectrum, and Clifford analysis \cite{Gab13, Ahn15, Vas15, Var15, Var16, Din17, Din17b, Pon17, Sle17,Vas17}
and they are therefore of special interest. With respect to higher spin
note that degree n-polynomials in $z\in \mathbb{C}$ form the carrier module for the $n+1$-dimensional irreducible spin representation of $SU(2,\mathbb{C})$ \cite{Man04}.

The relation to the previous sections can be drawn, if one considers polynomials in the solutions of the Laplace equation $\upsilon \in \mathbb{C}$.
If one operates on a polynomial in $\upsilon$ with $p_0=\partial_\upsilon$ and $q_0=\upsilon^2\partial_\upsilon$ the degree of the polynomial, and with this 
the spin of the irreducible representation of $SU(2,\mathbb{C})$, is modified.
In M\"obius and Lie sphere geometry there exists the notation $p_0\equiv e_\infty$ and $q_0\equiv e_0$, see for example Bobenko and Suris \cite{Bob07}. 
The reason is that $p_0$ maps on the projective line to infinity and $q_0$ to zero.

\section{Baryon number}
Houghton, Manton, and Sutcliffe use rational maps in the context of monopoles and the Skyrme model \cite{Hou98}. For skyrmions the degree of the rational map corresponds to the baryon number.
Meanwhile there have been calculations up to baryon number 108 \cite{Fei13}. The approach can be used also to calculate carbon states \cite{Lau14}.
Further applications of the rational map ansatz can be found in \cite{Bat14,Fos15}.
In this context the operators $q_0$ and $p_0$ and their complexifications can be used to modify the baryon number.

\section{Relativistic angular momentum tensor}
The solutions of the Laplace equation in $\mathbb{R}^{2}$ for $\alpha=1$ can be combined into two functions, which are defined 
using relations known from the definition of the sine and the cosine function
\be
\cn{\upsilon}=\frac{1}{2}\left(\upsilon+\frac{1}{\upsilon}\right),\quad \sn{\upsilon}=\frac{1}{2i}\left(\upsilon-\frac{1}{\upsilon}\right).
\ee
The motivation to introduce this generalization is given by the possibility
to represent the angular momentum tensor defined in Eq.~(\ref{coprot}) with the help of Eq.~(\ref{difform}) in a compact form
\be
\label{angularmomentum2}
s_{\mu\nu}=\left(\begin{array}{cccc}
0&i&i\sn{\upsilon}&-\cn{\upsilon}\\
-i&0&-i\cn{\upsilon}&-\sn{\upsilon}\\
-i\sn{\upsilon}&i\cn{\upsilon}&0&1\\
\cn{\upsilon}&\sn{\upsilon}&-1&0\\
\end{array}\right)\upsilon\partial_\upsilon\,.
\ee
Note that the tangent vector $\upsilon\partial_\upsilon$ is assigned to all elements of the Lie algebra.

The generalized cosine function is known as the Joukowski (Zhukovskii) map, which represents the complex potential function of a flow around a disk \cite{Olv15}. The disk can be used to model a quantum object
in the sense that matter represents a hole in the continuum of geometry. One may think here of the nucleon bag models \cite{Gui91}.
The generalized cosine function can thus be a useful mathematical tool, if one tries to describe quantum systems solely by geometry and topological invariants \cite{Man04,Wit88,Ati88}.

\section{The Hopf bundle}
The solutions of the Laplace equation can be applied also in the context of the Hopf map.
Because $\upsilon$ is just a complex number, one can use two solutions of the Laplace equation in the projective context
to define four real coordinate functions
\be
\upsilon_1=\varsigma_1+i\varsigma_2\,,\quad \upsilon_2=\varsigma_3+i\varsigma_4\,.
\ee
The two solutions of the Laplace equation in $\mathbb{R}^{2}$ can be identified with coordinates in $\mathbb{C}^2$ and thus also as a vector in Euclidean space $\mathbb{R}^{4}$
\be
\label{projrep}
\varsigma=\left(\begin{array}{c}
\upsilon_1\\
\upsilon_2\\
\end{array}\right).
\ee
The squared length 
of the vector can be expressed in terms of the real coordinates as
\be
\label{bilength}
\varsigma_1 \varsigma_1+\varsigma_2 \varsigma_2+\varsigma_3 \varsigma_3+\varsigma_4 \varsigma_4=\vert\varsigma\vert^2\,.
\ee
In the projective context the length $\vert \varsigma\vert$ has no relevance.
One can identify this vector then as a point on the sphere $S^3$. In this sense one can represent the Hopf map
$\pi: S^3\rightarrow S^2$ as \cite{Nak03}
\be
\label{Hopf}
\pi(\varsigma)=
\left(\begin{array}{c}
2(\varsigma_1 \varsigma_3+ \varsigma_2 \varsigma_4)\\
2(\varsigma_2 \varsigma_3- \varsigma_1 \varsigma_4)\\
\varsigma^2_1+ \varsigma^2_2- \varsigma^2_3- \varsigma^2_4
\end{array}\right)=\left(\begin{array}{c}
\xi_1\\
\xi_2\\
\xi_3
\end{array}\right)\,.
\ee
This transformation results in a vector $\xi$ which has the same length $\vert\xi\vert=\vert\varsigma\vert$ in the three-dimensional Euclidean metric $\mathbb{R}^{3}$.
Again the length has no relevance and the vector $\xi$ is considered as an element of $S^2$.

The Hopf bundle is defined in terms of two coordinate charts, which can be represented again by the two solutions of the Laplace equation in $\mathbb{R}^{2}$, now placed within the complex projective line $\mathbb{CP}^1$.
The  two coordinate charts, which relate the linear coordinates in $\mathbb{C}^2$ with the affine (stereographic) coordinates in $\mathbb{CP}^1$ are thus given by \cite{Mas96}
\be
\left(\begin{array}{c}
\upsilon_1\\
\upsilon_2\\
\end{array}\right)=\left(\begin{array}{c}
\upsilon_1/\upsilon_2\\
1
\end{array}\right),\quad
\left(\begin{array}{c}
\upsilon_1\\
\upsilon_2\\
\end{array}\right)=\left(\begin{array}{c}
1\\
\upsilon_2/\upsilon_1\\
\end{array}\right).
\ee
The details how the $U(1,\mathbb{C})=S^1$ transition functions between these two representations are derived can be found for example in \cite{Fra04}.
Due to the correspondence between $S^3$ and $SU(2,\mathbb{C})$ the Hopf bundle can be considered also as the principle bundle with total space $SU(2,\mathbb{C})$, and the fibration $SU(2,\mathbb{C})/U(1,\mathbb{C})=S^2$.

\section{Bicomplex numbers}
The research on bicomplex and hyperbolic numbers has been revived some decades ago by Yaglom \cite{Yag79}, Ryan \cite{Rya82}, and Price \cite{Pri91}.
Since then many applications of bicomplex and hyperbolic numbers have been published in mathematics and physics. 
Only a few of them are listed here \cite{Sob95, Roc04, Roc04b, Roc06, Cat08, Poo09, Lun12, Wan13, Alp14, Kum16, Gar16, Ban17, Ban17b, Kis17, Zar17}. Further
references can be found in these articles.

The two commutative complex units which form the bicomplex numbers are both identified with $\sqrt{-1}$.
In order to distinguish the two complex units, denoted here as $i$ and $j$, bicomplex generalizations of the null plane numbers \cite{Zho85,Huc93} can be introduced
\be
\label{derive}
i = o-\bar{o}\,,\quad j=o+\bar{o}\,.
\ee
This has the consequence that $i$ changes sign with respect to conjugation by means of the rule $\bar{\bar{o}}=o$ whereas $j$ remains invariant.
The calculation rules for the bicomplex null plane numbers can be derived from these definitions. The first set of rules has been introduced already in \cite{Ulr14}
\be
oo=i o\,,\quad \bar{o}\bar{o}=-i\bar{o}\,,\quad o\bar{o}=0\,.
\ee
There is a second complementary set of calculation rules, which can be used in this context
\be
oo=j o\,,\quad \bar{o}\bar{o}=j\bar{o}\,,\quad o\bar{o}=0\,.
\ee
One can consider two of the four square roots of the identity element. The square roots result by multiplication of the above bicomplex units, e.g., $ii=\sqrt{(-1)(-1)}=\sqrt{1}$
\be
ij=\sqrt{1}\,,\quad -1=\sqrt{1}\,.
\ee
The first of the square roots $ij$ can be identified with the hyperbolic unit, which can be used to construct the hyperbolic numbers.
The square roots of the identity can be represented again in terms of the bicomplex null plane numbers
\be
 ij=oo-\bar{o}\bar{o}\,,\quad -1=oo+\bar{o}\bar{o}\,.
\ee
One finds that the hyperbolic unit $ij$ changes sign with respect to conjugation, whereas the negative identity does not change its sign as expected.

\section{Bicomplex numbers and the Hopf bundle}
\label{bihopf}
With the help of the hypercomplex units described in the previous section the bicomplex number
can be defined as
\be
\varsigma=\varsigma_1+i\varsigma_2+j\varsigma_3+ij\varsigma_4\,.
\ee
The squared length of the bicomplex number will be defined again by Eq.~(\ref{bilength}).
In the projective context the length itself is an irrelevant quantity. Thus the bicomplex number can be interpreted as a point on the sphere $S^3$.
With conjugation defined in the preceding section one can calculate in detail
\be
\label{projconj}
\varsigma\bar{\varsigma}=\xi_3+j\xi_1\,.
\ee
The map $\varsigma\bar{\varsigma}$ is not leading to the squared length. 
Conjugation projects instead towards two coordinates in the base space of the Hopf bundle, see Eq.~(\ref{Hopf}). 

The third coordinate $\xi_2$  can be determined by a second involution, the reversion
\be
\varsigma^\dagger=\varsigma_1-i\varsigma_2-j\varsigma_3+ij\varsigma_4\,.
\ee
The sign change of the hypercomplex units is consistent with previous publications of the author on hyperbolic numbers \cite{Ulr05}.
The reader should keep in mind that in contrast to \cite{Ulr05} the hyperbolic unit is now denoted as $ij$ instead of $j$.
One can calculate the following relationship
\be
\varsigma\varsigma^\dagger=\vert\varsigma\vert^2-ij\xi_2\,.
\ee
Thus conjugation and reversion can be
used to draw a relation to the Hopf map and the squared length of the bicomplex number.
It is therefore also possible to identify the bicomplex number with an element of the complex projective line
\be
\label{projcont}
\varsigma=\upsilon_1+j\upsilon_2\,.
\ee
The notation indicates how the bicomplex numbers are mapped to the representation defined in Eq.~(\ref{projrep}).

Up to scale the bicomplex numbers correspond to the sphere $S^3$ like the complex numbers correspond to $S^1$.
From a physical point of view one can place the mathematics of the Dirac monopole into this context because the total space of the corresponding non-trivial fiber bundle corresponds to $S^3$ \cite{Goc87}.
In the base space the pole structure of the monopole potential refers to the matrix $q_0$ as it appears in Eq.~(\ref{qpole}).
In this sense one can use the Euclidean space $\mathbb{R}^4$ also to consider instantons. As discussed the space is made up of generalized coordinates $\varsigma$ related to the solutions $\upsilon$ of the Laplace equation.

\section{Bicomplex spin representation}
\label{bispin}
The bicomplex numbers as discussed in the previous sections can be used
to generalize the representation of the matrices in Eq.~(\ref{realproj}) to the bicomplex $S^3$ context. This is necessary to describe relativistic quantum systems as will be discussed below.

The new representation of $p_0$ and $q_0$ is constructed with linear combinations of the original matrices and multiplication by the bicomplex units $i$ and $j$
\be
p_0=\left(\begin{array}{cc}
0&o\\
-\bar{o}&0
\end{array}\right),\quad
q_0=\left(\begin{array}{cc}
0&\bar{o}\\
-o&0
\end{array}\right).
\ee
With the help of the relevant Lie algebra relation in Eq.~(\ref{confcomm}) one can derive the
representation of the scale operator, which corresponds to a multiplication of the original matrix by the hyperbolic unit $ij$
\be
\label{hypscale}
b=\frac{1}{2}\left(\begin{array}{cc}
ij&0\\
0&-ij
\end{array}\right).
\ee
The remaining matrices of the relativistic spin operator can be derived from Eq.~(\ref{complexline}).
The result is consistent with \cite{Ulr17}.

Bicomplex and hyperbolic numbers 
allow for a more compact notation if relativistic negative energy quantum states are considered \cite{Ulr05}.
In this sense the bicomplex $2\times 2$ matrices coincide with the dimensionality of the complex $4\times 4$ Dirac matrices which take into account discrete transformations like the charge, parity, and time transformations \cite{Huc93}.

\section{Clifford paravector algebra in terms of tangent vectors}
From the spin representations one can return now to the two-dimensional plane, where the Laplace equation has been initially defined.
The complex numbers can be understood in this context as referring to the
Clifford paravector algebra $\mathbb{R}_{0,1}$ with the basis elements $e_\mu=(e_0,e_1)=(1,i)$. Clifford algebras are discussed in detail for example in \cite{Por81,Por95,Lou02,Vaz16}.
One can identify the $1\times 1$ matrix representation of the one parameter transformation generated by $e_0=1$ with the scale transformation
\be
\label{nscale1}
e^{e_0\varepsilon}=e^\varepsilon\,.
\ee
A point on the group manifold of dilations and
rotations can be represented by the solution of the Laplace equation $\upsilon$.
As in Exercise 9.13 of Misner, Thorne, and Wheeler \cite{Mis73} one can consider a curve $c(\varepsilon)$ on this group manifold
\be
\label{nscale2}
c(\varepsilon)=e^\varepsilon\upsilon\,.
\ee
To give an approximation for this curve the exponential $e^\varepsilon\approx 1+\varepsilon$ is applied to the Taylor expansion of the sine function included in $\upsilon$.
The sine function can be expanded around an arbitrary angle $\theta_0$. One finds the result
\be
c(\varepsilon)\approx\sin({\theta+\varepsilon\tan{\theta}})e^{i\varphi}\,.
\ee
The tangent vector $e_0$ at the point $\upsilon$ of the group manifold is then given by
\be
e_0\upsilon=\frac{d}{d\varepsilon}c(\varepsilon)\vert_{\varepsilon=0}=\tan{\theta}\partial_\theta \upsilon\,.
\ee
In a similar way one finds for the tangent vector along the curve generated by the $1\times 1$ matrix representation of the second basis element $e^{e_1\varepsilon}=e^{i\varepsilon}$
\be
e_1\upsilon=\frac{d}{d\varepsilon}e^{i\varepsilon}\upsilon\vert_{\varepsilon=0}=\partial_\varphi \upsilon\,.
\ee
The elements of the Clifford paravector algebra can thus be identified with their representation in terms of
tangent vectors $e_\mu=(e_0,e_1)$ as 
\be
\label{equirep}
e_\mu=(1,i)=(\tan{\theta}\partial_\theta,\partial_\varphi)=(b,s_{01})\,. 
\ee
The next two sections will make use of this identification to obtain further insights with respect to the solutions of the Laplace equation.

\section{Fock space}
\label{Fock}
The solution of the Laplace equation can be reformulated using the correspondence derived in the previous section.
In a first step Eq.~(\ref{SolofLap}) can be separated for $\alpha=1$ into real and complex contributions
\be
\upsilon=\sin{\varphi}\sin{\theta}+i\sin{\varphi}\sin{\theta}\,.
\ee
The solution of the Laplace equation is now interpreted as a Clifford paravector expressed in the basis $e_\mu$.
By means of Eq.~(\ref{equirep}) it is now possible to replace the identity in the first term and the complex unit in the second term, which results in
\be
\upsilon=\sin{\varphi}\sin{\theta}\tan{\theta}\partial_{\theta}+\sin{\varphi}\sin{\theta}\partial_\varphi=q_0.
\ee
The solution of the Laplace equation can be identified directly with $q_0$, one of the generators of the symmetry group of the Laplace equation as given by Eq.~(\ref{allop}). 
In the same way Eq.~(\ref{equirep}) can be used
to calculate
\be
\upsilon^{-1}=p_0\,.
\ee
It is therefore possible to identify $p_0$ with the inverse of the solution.

One can compare these results with the corresponding representations on the projective line $\mathbb{CP}^1$, which are provided by Eq.~(\ref{difform}) 
\be
\label{vecfund}
q_0=\upsilon (\upsilon\partial_\upsilon)\,,\quad p_0=\upsilon^{-1}(\upsilon\partial_\upsilon)\,.
\ee
Thus one finds that $q_0$, which maps an element of the projective line to 0, can be interpreted as a creation operator. In contrast $p_0$ as the map to infinity annihilates one unit of the solution of the Laplace equation.

\section{Cartan geometry}
\label{cartgeo}
In \cite{Ulr17} the M\"obius geometry is defined based on translations as the coset representatives of the corresponding homogeneous space.
Dilations and rotations represent the commuting Cartan subalgebra of the conformal algebra, which can also be used to represent the Euclidean subgeometry as discussed in Section \ref{geointeract}.
This suggests to change the homogenous space accordingly.
One can work now with
\be
\label{coset}
\mathfrak{g}=\{e_\mu, p_\mu,q_\mu\}\,,\quad\mathfrak{h}=\{p_\mu,q_\mu\}\,.
\ee
Again one should note that $e_\mu=(e_0,e_1)=(b,s_{01})$. 
Dilations and rotations are now the coset representatives of a Klein geometry, which can be further extended to a Cartan geometry as discussed by Sharpe \cite{Sha97}.

The method can be further developed using differential forms which are dual to the vector fields considered so far.
To provide a compact notation one can look at the representation of the spin operator in Eq.~(\ref{angularmomentum2})
and require for a dual differential form $\omega$
\be
\omega^{\mu\nu}(s_{\rho\sigma})=\delta^\mu_{\phantom{\mu}\rho}\delta^\nu_{\phantom{\nu}\sigma}\,.
\ee
Here it is referred to the notation used in \cite{Goc87}. The differential form $\omega$ allows to further proceed with the development of a theory that extends general relativity following the methods discussed by Sharpe \cite{Sha97}.
The higher dimensional ambient space $\mathbb{R}^{3,1}$ appears via the spin angular momentum and Eq.~(\ref{comm}), which allows to study aspects of the AdS/CFT correspondence in a simplified geometry.

\section{Rotations and dilations on the projective line}
So far the discussion restricts to a pre-gauge theory. It is beyond the scope of this article to develop the gauge theory in detail as discussed in the previous section.
Nevertheless, to motivate investigations into this direction an attempt is made to relate the coset representatives $e_\mu$ of the coset space $G/H$, defined by Eq.~(\ref{coset}) in the context of general relativity, to generators of gauge groups in 
the corresponding conformal field theory.

One can consider the coset representatives $e_\mu$ acting with their spin representations on the complex projective line.
The explicit representation of $s_{01}$
can be found in \cite{Ulr17}. One can exponentiate this expression and then operate with $\exp{(s_{01}\varepsilon)}$ on a projective vector. In detail one finds
\be
\label{rotsym}
\left(\begin{array}{cc}
e^{i\varepsilon/2}&0\\
0&e^{-i\varepsilon/2}
\end{array}\right)
\left(\begin{array}{c}
\upsilon_1\\
\upsilon_2\\
\end{array}\right)=
\left(\begin{array}{c}
e^{i\varepsilon}\upsilon_1/\upsilon_2\\
1\\
\end{array}\right)=\left(\begin{array}{c}
e^{i\varepsilon}\upsilon\\
1\\
\end{array}\right).
\ee
Here the linear coordinates in $\mathbb{C}^2$ on the left hand side are changed into the stereographic coordinates of $\mathbb{CP}^1$ on the right hand side.
The first equality includes the rescaling to set the lower component to the identity. In the second equality the redefinition $\upsilon=\upsilon_1/\upsilon_2$ is used.
In the same way the spin representation for $b$ as derived in Eq.~(\ref{hypscale}) can be inserted and the exponential $\exp{(b\varepsilon)}$ is applied to the projective vector
\be
\label{baseexponent2}
\left(\begin{array}{cc}
e^{ij\varepsilon/2}&0\\
0&e^{-ij\varepsilon/2}
\end{array}\right)
\left(\begin{array}{c}
\upsilon_1\\
\upsilon_2\\
\end{array}\right)=\left(\begin{array}{c}
e^{ij\varepsilon}\upsilon\\
1\\
\end{array}\right).
\ee
Instead of using the bicomplex representation of the scale operator $b$ one can insert the standard representation given by Eq.~(\ref{realproj})
into the exponential. This leads via Eq.~(\ref{baseexponent2}) to a representation of the scale operator which is consistent with the familiar expression
as used in Eqs.~(\ref{nscale1}) and (\ref{nscale2}).

\section{Maxwell theory of gravitation}
The two hypercomplex units visible in the exponents of rotations and dilations on the right hand side of Eqs.~(\ref{rotsym}) and (\ref{baseexponent2}) can now be used to set up two gauge symmetries.
The gauge theory related to $i$ is identified with electromagnetism. The generator $ij=\sqrt{1}$ can be related to
a Maxwell theory of gravitation as discussed in \cite{Ulr06}.
The perihelion 
of a planet, the deflection of light in the gravitational field of a star, and the 
gravitational red shift, as predicted by the results of the general theory of relativity, have been calculated in a Maxwell theory of gravitation for example by Majern\'ik \cite{Maj72,Maj82} and Singh \cite{Sin82}.
Singh achieves this by using a four potential $\phi_\mu$, which includes the Newton potential as its time component $\phi_0\sim 1/r$ and additional velocity dependent contributions $\phi_i\sim v_i/r$
in the vector part of the potential.

\section{Conformal coordinates}
\label{confrep}
As dilations and rotations have now been identified to be of central importance in the context of the Laplace equation
it is worth to look at coordinates, which puts the corresponding generators in their most simple form
\be
\label{cpara}
x_{\mu}=\left(\begin{array}{c}
x_{0}\\
x_{1}\\
\end{array}\right)=
\left(\begin{array}{c}
e^\rho\cos{\varphi}\\
e^\rho\sin{\varphi}\\
\end{array}\right).
\ee
Here $\rho$ parameterizes the scale variable and $\varphi$ the usual angle coordinate in the two-dimensional plane.
The resulting basis vectors can be used to calculate the metric tensor
\be
g^{\alpha\beta}= \varsigma^\alpha \cdot \varsigma^\beta=\left(\begin{array}{cc}
e^{2\rho}&0\\
0&e^{2\rho}\\
\end{array}\right).
\ee
In the context of conformal geometry all these metrics are considered to be equivalent, see for example Folland \cite{Fol70}.
The conformal coordinates thus refer to the ideas of Weyl \cite{Wey18}. However, the physical interpretation is now different. The two-dimensional space as an Euclidean worldsheet corresponds to a geometric polarization plane,
which resides in arbitrarily rotated form within Minkowski space \cite{Ulr17}. The Euclidean transformations within the polarization plane, represented by rotations and dilations, 
give rise in this model to electromagnetism and gravitation as discussed in the previous two sections.

The transformation matrix referring to the conformal coordinate system can be calculated as defined in Eq.~(\ref{cindextrafo})
\be
\label{lamtrafo}
A^\mu_{\phantom{\mu}\alpha}=\left(\begin{array}{cc}
e^{-\rho}\cos{\varphi}&-e^{-\rho}\sin{\varphi}\\
e^{-\rho}\sin{\varphi}&e^{-\rho}\cos{\varphi}\\
\end{array}\right).
\ee
As outlined in the preceding sections this operator can be used to calculate representations of the generators in the new coordinates.
The coordinates introduced above result in the simplest representation of scale and rotation operators
\be
b=\partial_\rho\,,\quad s_{01}=\partial_{\varphi}\,.
\ee
This representation can be used to formulate the Laplace equation,
which is also found to be in the simplest form beside the standard cartesian coordinate representation
\be
\label{conlap}
\triangle=\partial_{\rho\rho}
+\partial_{\varphi\varphi}\,.
\ee
The solutions coincide again with the solutions of the representation in cartesian coordinates.
According to the discussion in Sections \ref{bihopf} and \ref{bispin}, the solutions can again be considered either in the complex $S^1$ or the relativistic bicomplex $S^3$ context 
\be
\upsilon^\alpha=e^{i\alpha\varphi}e^{\alpha\rho}\,,\quad \upsilon^\alpha=e^{i\alpha\varphi}e^{ij\alpha\rho} \,. 
\ee
The conformal coordinates are thus better adjusted to the impact of rotations and dilations. Compared to the holographic coordinates 
this representation suffers however from the non-compactness of the parameter space with regards to~$\rho$.

\section{Summary}
The Laplace equation in the two-dimensional Euclidean space serves as an example on how a coordinate parametrization can look like, if it does not
rely on length scales. The representation is based on the holographic coordinates, which refer to the stereographic projection and the conformal transformations.
The solutions of the Laplace equation in holographic coordinates can be identified
with a subset of the spherical harmonics in three-dimensional Euclidean space, which connects the method with higher dimensional geometries.

Concepts in geometry and physics like projective spaces, the Hopf map, relativistic spin, bicomplex numbers, and the Fock space can
be connected through the Laplace equation by means of complex coordinates made up of the solutions of the Laplace equation. 
The complex numbers as the most simple Clifford paravector algebra are reinterpreted as tangent vectors.
Using this interpretation the relation to Cartan geometries is outlined, which can serve as the starting point to investigate the gauge theory of the method.
Here dilations and rotations play the central role.

\end{document}